# A meshless numerical solution of the family of generalized fifth-order Korteweg-de Vries equations

Syed Tauseef Mohyud-Din[1], Elham Negahdary[2], Muhammad Usman[2]
[1]HITEC University, Taxila Cantt, Pakistan,
[2]University of Dayton,
300 College Park,
Dayton, OH 45469-2316, USA.

**Abstract**
In this paper we present a numerical solution of a family of generalized fifth-order Korteweg-de Vries equations using a meshless method of lines. This method uses radial basis functions for spatial derivatives and Runge-Kutta method as a time integrator. This method exhibits high accuracy as seen from the comparison with the exact solutions.

## 1. Introduction

Most of the physical phenomena in nature are modeled by nonlinear partial differential equations (PDEs). In general, nonlinear PDEs cannot be solved analytically and hence need to be solved numerically in order to predict the behavior of the system.

Commonly used numerical methods to approximate the solutions of nonlinear PDEs include finite difference methods, collocation methods and Galerkin methods. However, some of these methods are not easy to use and sometimes require tedious work and calculation [32, 3].

Finite-difference methods are known as effective tools for solving a variety of PDEs [6]. Conditional stability of explicit finite-difference schemes puts a severe constraint on the time step, while implicit finite-difference schemes are computationally expensive [4]. Furthermore, these methods can be made highly accurate, but require a structured grid. Another class of methods known as spectral methods is found to be even more accurate but have restriction on the geometry of the problem, for example in the Fourier case we need to have periodic boundary conditions. Finite-element methods have been used as an alternative method for numerical solution of PDEs. This family of numerical techniques is efficient particularly for solving problems with arbitrary geometry. But the need to produce a body-fitted mesh in two- and three-dimensional problems makes these methods quite time-consuming and difficult to use [1]. Overall, finite-element techniques are highly flexible, but it is hard to obtain results with high-order accuracy.

Consequently, to avoid the mesh generation, meshless techniques are attracted the attention of researchers, in recent years, as alternatives to traditional finite element, finite volume and finite difference methods. In a meshless (meshfree) method a set of scattered nodes, with no connectivity information required among the set of points, is used instead of meshing the domain of the problem. Examples of some meshless schemes are the element free Galerkin method, the reproducing kernel particle, the local point interpolation, etc (e.g. see [23] and references therein).

Over the last two decades, the radial basis function methods have emerged as a powerful tool for scattered data interpolation problems. The use of radial basis functions as a meshless procedure for numerical solution of PDEs is based on the collocation scheme. Due to the collocation technique, this method does not need to evaluate any integral. The main advantage of numerical procedures that use radial basis functions over traditional techniques is the meshless property of these methods. Radial basis functions are actively used for solving PDEs ([5, 18], and references therein). In the above cited work, RBFs are used to replace the function and its spatial derivatives, while a finite difference scheme is used to march in time.

This method was first introduced by Kansa [19, 20] in 1990 for the numerical solutions of the PDEs. Kansa used the Multiquadric (MQ) RBF to solve the elliptic and parabolic PDEs. Recently, Flyer and Wright [9] indicated RBFs allowed for a much lower spatial resolution, while being able to take unusually large time-steps to achieve the same accuracy compared to other methods.

In this paper we will use method of lines coupled with RBFs to find the numerical solution of the family of generalized fifth order KdV (henceforth gfKdV) equation:

$$u_t + au^2u_x + bu_xu_{xx} + cuu_{xxx} + du_{xxxxx} = 0, \tag{1}$$

where $a$, $b$, $c$ and $d$ are constants. Eq. (1) is known as Lax's fifth-order KdV equation [31], if we set $a$ = 30, $b$ = 30, $c$ = 10 and $d$ = 1 and the Sawada-Kotera equation for $a$ = 45, $b$ =15, $c$ = 15, $d$ = 1 [25, 21].

The method of lines (MOL) [29] is generally recognized as a comprehensive and powerful approach to the numerical solution of time-dependent PDEs. This method is comprised of two steps: first, approximating the spatial derivatives; Second, then resulting system of semi- discrete ordinary differential equations (ODEs) is integrated in time. Hence the method of lines approximates the solution of PDEs using ODEs integrators.

In this paper, we will use radial basis functions combined with the MOL, hence calling it MOL-RBF, to solve the gfKdV equation inspired by [30]. As evident from our results, this method possesses high accuracy and ease of implementation. The computed results are compared with the analytic solutions and good agreement is indicated. The remainder of the paper is organized as follows. In Section 2, we show the formulation of RBF method and then, we couple the RBFs meshless method with the MOL to solve the gfKdV equations. In Section 3, we apply this method to Lax's fifth-order KdV equation and the Sawada-Kotera equation as two examples of gfKdV. Section 4 is a comparison of our results with the exact solution and analysis of our method. The last section is a brief conclusion.

## 2. The MOL-RBF method

A radial basis function is a kind of function with the independent variable $r_i = r(x, x_i) = \| x - x_i \|$. Some of the commonly used RBFs in the literature are:

$\phi_i(x) = (c^2 + r_i^{\,2})^{1/2}$ Multiquadric (MQ),

$\phi_i(x) = (c^2 + r_i^{\,2})^{-1/2}$ Inverse Multiquadric (IMQ),

$\phi_i(x) = e^{-cr_i^2}$ Gaussian (GA),

where the free parameter $c$ is called the shape parameter of the RBFs. In the above definition x = ($x$, $y$) are the cartesian coordinates in $\Omega \subset R^2$ and the radius is given by

$$r_j = \| \mathbf{x} - x_j \| = \left\{ (x - x_j)^2 + (y - y_j)^2 \right\}^{\frac{1}{2}},$$

where $(x_j, y_j)$ is called the $j^{th}$ source point of the RBF and is denoted by $\mathbf{x}_j$. We choose N nodes $(x_i,\ i = 1,2,...,N)$ in $\wp_s \subseteq \Omega \cup \partial\Omega$. Any given smooth function can be represented as a linear combination of RBFs:

$$u^N(x_i) = \sum_{i=1}^{N} \lambda_i \phi_i = \mathbf{\Phi}^T(x)\lambda, \tag{2}$$

where

$$\mathbf{\Phi}(x) = [\phi_1(x), \phi_2(x), ..., \phi_N(x)]^T,$$

$$\boldsymbol{\lambda} = [\lambda_1, \lambda_2, ..., \lambda_N]^T.$$

Which can be written as

$$\mathrm{A}\lambda = \mathrm{u},$$

where

$$\mathbf{u} = [u_1, u_2, ..., u_N]^T,$$

and the matrix

$$\mathbf{A}=\begin{bmatrix}\mathbf{\Phi}^T(x_1)\\ \mathbf{\Phi}^T(x_2)\\ \cdots\\ \mathbf{\Phi}^T(x_N)\end{bmatrix}=\begin{bmatrix}\phi_1(x_1) & \phi_2(x_1) & \cdots & \phi_N(x_1)\\ \phi_1(x_2) & \phi_2(x_2) & \cdots & \phi_N(x_2)\\ \cdots & \cdots & \ddots & \cdots\\ \phi_1(x_N) & \phi_2(x_N) & \cdots & \phi_N(x_N)\end{bmatrix}$$

is called the interpolation matrix, consisting of functions forming the basis of the approximation space. It follows from Eq. (2) and $A\lambda = u$, that

$$u^N(x) = \mathbf{\Phi}^T(x)\mathbf{A}^{-1}\mathbf{u} = \mathbf{V}(x)\mathbf{u},$$

where

$$\mathbf{V}(x) = \mathbf{\Phi}^T(x)\mathbf{A}^{-1}.$$

The convergence of RBF interpolation is given by the theorems in [33, 34]:
Assuming $\{x_i\}_{i=1}^N$ are $N$ source points in $\wp_s$ which is convex, the radial distance is defined as

$$\delta := \delta(\Omega, \wp_s) = \max_{x\in\Omega}\min_{1\le i\le N}\|x - x_i\|_2, \tag{3}$$

we have

$$\|u^N(x) - u(x)\| \approx O(\eta^{c/\delta}), \tag{4}$$

where $0 < \eta < 1$ is a real number and $\eta = \exp(-\theta)$ with $\theta > 0$.
From (4) it is clear that the parameter $c$ and radial distance $\delta$ the rate of convergence.
The exponential convergence proofs in applying RBFs in Sobolov space was given by Yoon [35], spectral convergence of the method in the limit of flat RBFs was given by Fornberg et al. [10]. The exponential convergence rate was verified numerically by Fedseyev et al. [8]. The exponential convergence cited above is limited to certain classes of functions that are smooth enough and well-behaved in the domain of approximation.
In 1971, Hardy [15] developed multi-quadric MQ to approximate two-dimensional geographical surfaces. In Franke's [12] review paper, the MQ was rated one of the best methods among 29 scattered data interpolation schemes based on their accuracy, stability, efficiency, ease of implementation, and memory requirement. Further, the interpolation matrix for MQ is invertible. In 1990, since Kansa [19, 20] modified the MQ for the solution of elliptic, parabolic and hyperbolic type PDEs, radial basis functions has been used to solve partial differential equations numerically [2, 17, 18, 26, 30]. The accuracy of MQ depends on the choice of a user defined parameter $c$ called the shape parameter that affects the shape of the RBFs. Golberg, Chen, and Karur [13] and Hickernell and Hon [16] applied the technique of cross validation to obtain an optimal value of the shape parameter c.
The non-singularity of the collocation matrix **A** depends on the properties of RBFs used. According to [24], the matrix **A** is conditionally positive definite for MQ radial basis functions. This fact guarantees the non-singularity of the matrix A for distinct supporting points.
Now, we apply the method of lines combined with the RBFs (MOL-RBFs) for gfKdV.

$$u_t + au^2u_x + bu_xu_{xx} + cuu_{xxx} + du_{xxxxx} = 0, \quad x\in[a,b]$$

with the following initial condition and boundary conditions

$$u(x,t_0) = u_0(x), u(a,t) = f(t), u_x(a,t) = 0, u(b,t) = u_x(b,t) = u_{xx}(b,t) = g(t),$$

where a, b, c and d are real constants and $u_0(x)$, $f(x)$ and $g(x)$ are known functions.
First, we choose $N$ nodes in [a, b]

$$a = x_1 < x_2 < ... < x_{N-1} < x_N = b.$$

By RBF interpolation, we get

$$u(x,t) \approx u^N(x,t) = \sum_{i=1}^N \lambda_i\phi_i = \mathbf{\Phi}^T(x)\mathbf{A}^{-1}u = \mathbf{V}(x)\mathbf{u},$$

where

$$\mathbf{\Phi}(x) = [\phi_1(x), \phi_2(x), ..., \phi_N(x)]^T,$$

$\mathbf{u} = [u_1(t), u_2(t), ..., u_N(t)]^T$,

$\mathbf{V}(x) = \mathbf{\Phi}^T(x)\mathbf{A^{-1}} = [V_1(x)\ V_2(x), ..., V_{N-1}(x)\ \ V_N(x)]$

By applying this method to the gfKdV, we obtain

$$\frac{du_i}{dt} + au_i^2(\mathbf{V}_x(x_i)\mathbf{u}) + b(\mathbf{V}_x(x_i)\mathbf{u})(\mathbf{V}_{xx}(x_i)\mathbf{u}) + cu_i(\mathbf{V}_{xxx}(x_i)\mathbf{u}) + d(\mathbf{V}_{xxxxx}(x_i)\mathbf{u}) = 0,\ \ i = 1, 2, ..., N, \tag{5}$$

where $u_i$ is abbreviation of $u_i(t)$

$\mathbf{V_x}(x_i) = [V_{1x}(x_i)\ V_{2x}(x_i), ..., V_{Nx}(x_i)]$,

$V_{jx}(x_i) = \frac{\partial}{\partial x}(V_j(x_i))$,

$\mathbf{V_{xx}}(x_i) = [V_{1xx}(x_i)\ V_{2xx}(x_i), ..., V_{Nxx}(x_i)]$,

$V_{jxx}(x_i) = \frac{\partial^2}{\partial x^2}(V_j(x_i))$,

$\mathbf{V_{xxx}}(x_i) = [V_{1xxx}(x_i)\ V_{2xxx}(x_i), ..., V_{Nxxx}(x_i)]$,

$V_{jxxx}(x_i) = \frac{\partial^3}{\partial x^3}(V_j(x_i))$,

$\mathbf{V_{xxxxx}}(x_i) = [V_{1xxxxx}(x_i)\ V_{2xxxxx}(x_i), ..., V_{Nxxxxx}(x_i)]$,

$V_{jxxxxx}(x_i) = \frac{\partial^4}{\partial x^4}(V_j(x_i))$,

If we apply the collocation to gfKdV, Eq. (1) will take the form

$$\frac{dU}{dt} + a\mathbf{U^2} * (M_x\mathbf{U}) + b(M_x\mathbf{U}) * (M_{xx}\mathbf{U}) + c\mathbf{U} * (M_{xxx}\mathbf{U}) + d(M_{xxxxx}\mathbf{U}) = 0 \tag{6}$$

where $\mathbf{U} = [u_1\ u_2, ..., u_{N-1}\ u_N(t)]^T$,

$M_x = [V_{jx}(x_i)]_{N\times N}$,

$M_{xx} = [V_{jxx}(x_i)]_{N\times N}$,

$M_{xxx} = [V_{jxxx}(x_i)]_{N\times N}$,

$M_{xxxxx} = [V_{jxxxxx}(x_i)]_{N\times N}$,

and $*$ denotes the component by component multiplication of two vectors. We can rewrite (6) as:

$H(\mathbf{U}) = -a\mathbf{U^2} * (M_x\mathbf{U}) - b(M_x\mathbf{U}) * (M_{xx}\mathbf{U}) - c\mathbf{U} * (M_{xxx}\mathbf{U}) - d(M_{xxxxx}\mathbf{U})$.

The initial condition vector is

$\mathbf{U}(t_0) = [u_0(x_1), u_0(x_2), ..., u_0(x_N)]^T$,

and we use two dirichlet boundary conditions;

$u_1(t) = f(t),\ \ u_N(t) = g(t)$.

It is obvious that our PDE has become an ODE and this ODE can be solved by any of several ODE solvers, we choose the fourth order Runge-Kutta scheme (RK4)

$$U^{n+1} = U^n + \frac{\Delta t(K_1 + 2K_2 + 2K_3 + K_4)}{6}$$

$K_1 = H(U^n)$,

$K_2 = H(U^n + \frac{\Delta t}{2}K_1)$,

$K_3 = H(U^n + \frac{\Delta t}{2}K_2)$,

$K_4 = H(U^n + \Delta t\, K_3)$,

Now we are ready to apply this method to solve gfKdV numerically.

## 3. Application

We apply the MOL-RBF method to Lax's and Swada-Kotera cases of gfKdV equations.

**Lax's fifth-order KdV equation**

$$u_t + 30u^2 u_x + 30u_x u_{xx} + 10uu_{xxx} + u_{xxxxx} = 0 \qquad (7)$$

with the initial condition

$u(x,0) = 2k^2(2 - 3\tanh^2(k(x - x_0)))$.

The exact solution is given by

$u(x,t) = 2k^2(2 - 3\tanh^2(k(x - 56k^4 t - x_0)))$.

**Sawada-Kotera (SK) equation**

Next, we consider the Sawada-Kotera equation

$$u_t + 45u^2 u_x + 15u_x u_{xx} + 15uu_{xxx} + u_{xxxxx} = 0 \qquad (8)$$

with initial condition is given by

$u(x,0) = 2k^2 \operatorname{sec}h^2(k(x - x_0))$.

The exact solution is given by

$u(x,t) = 2k^2 \operatorname{sec}h^2(k(x - 16k^4 t - x_0))$.

We apply the numerical method and evaluate and compare the solutions with exact solutions and present the results in Tables1- 4. We use the Max-error, $L2$-error and Root Mean Square (RMS) error

$$\|u^N - u\|_{\max} = \max_{1 \le i \le N} |u_i^N - u|,$$

$$\|u^N - u\|_{L^2} = \sqrt{h \sum_{i=1}^{N} (u_i^N - u)^2},$$

$$\text{RMS error} = \sqrt{\frac{\sum_{i=1}^{N} (u_i^N - u)^2}{N}},$$

where $u^N$ is the approximate solution and $u$ is the exact solution of Eq. (1)

## 4. The analysis of numerical experiment

We have applied the meshless method of lines (MOL-RBF) to Lax and SK. Figures 1 and 2 show the 3d plots of computed solutions of the gfKdV equation for Lax and SK cases respectively. We have omitted the 3d plots of exact solutions as they appear identical to the computed solutions. Figures 4- 7 show Lax and SK solutions at $t$ =2 for $k$ = 0.001 and $k$ = 0.00001, for different RBFs. Since the accuracy depends upon number of nodes and the value of the shape parameter $c$. Choosing the optimal value of shape parameter is still an open problem. Many researchers proposed methods to find the optimal $c$. For example Hardy's formula

$c = 0.815d$ where $d = \frac{1}{N}\sum_{i=1}^{N} d_i$

where $d_i$ is the distance from the $i^{th}$ center to the nearest neighbor, $N$ is the number of centers.

Franke's formula for the shape parameter is

$$c = \frac{1.25D}{\sqrt{N}},$$

where $D$ is the diameter of the smallest circle encompassing all the centers and as before $N$ represents the total number of centers.

According to [27] the shape parameter should depend on factors like: number of grid points, distribution of grid points, interpolation function $\phi$, condition number of the matrix, and computer precision. We have used the very common brute force method to find the optimal value of the shape parameter. In this method we plotted the error vs. shape parameter and then picked *c* corresponding to minimum error and not too high condition number. We have solved the gfKdV equation using the meshless method of lines for the values $x_0 = 0,\ k = 0.001$ and $x_0 = 0,\ k = 0.00001$. We have presented our results of Max-error, $L_2$-error and RMS error in Tables 1 and 2. For both cases (Lax and SK) the optimal value of the shape parameter is of the same order of magnitude for GA, MQ and IMQ.

One more point to be noted here is that the condition number is approximately the same for both Lax and SK in each of GA, MQ and IMQ case. Table 2 also shows a trade off for between the accuracy and condition number for the radial basis function IMQ, this is termed as *uncertainty* or *trade-off principle* in [11]. Figure 3 is the plot for the solution of gfKdV (Lax) using GA radial basis function and non-optimal value of shape parameter *c* = 3690. It is clear that the error is higher at the end of the interval as pointed in [11, 22] for non-optimal shape parameters. In Table 3 and 4 we have computed the conserved densities [14], $I_1 = \int u\,dx$ (Lax and SK), $I_2 = \int \frac{1}{3}u^3 - u_x^2\,dx$ (SK), $I_2 = \int \frac{1}{3}u^3 - \frac{1}{6}u_x^2\,dx$ for Lax.

## 5. Conclusion

In this paper, we have solved the family of generalized-fifth order Korteweg-de Vries (gfKdV) equations. We proposed a meshless MOL with the use of RBFs for solving the gfKdV equations. The numerical results given in the previous section demonstrate the good accuracy of this method. It has to be emphasized that the shape parameter for all the calculations performed in this paper was found experimentally. Results of numerical experiments indicate that the Gaussian RBF has the best accuracy in this method for gKdV equation with the condition number of exactly one. Also, the RBFs allow for a much lower spatial resolution (i.e., lower number of nodes) to obtain much higher accuracy. The two major advantages of this method are the meshless property and use of ODE solvers of high quality and their codes to approach to the solutions of PDEs. This method can be extended to solve the PDE with higher order derivatives with respect to *x*, e.g., Kuramoto-Sivashinsky equation, without any difficulties because the radial basis functions MQ, IMQ, and GA are infinitely differentiable. We consider this method as one of the more efficient methods and easy to implement for the numerical solutions of PDEs. We have used the fixed value of shape parameter that minimizes the error because of known exact solution. Interested readers can try a different approach for example optimization via residual error to choose the shape parameter if analytical solution of the problem is not known or recently [28] shows the use of random values for the shape parameter and in [7] adaptive residual subsampling method has been proposed.

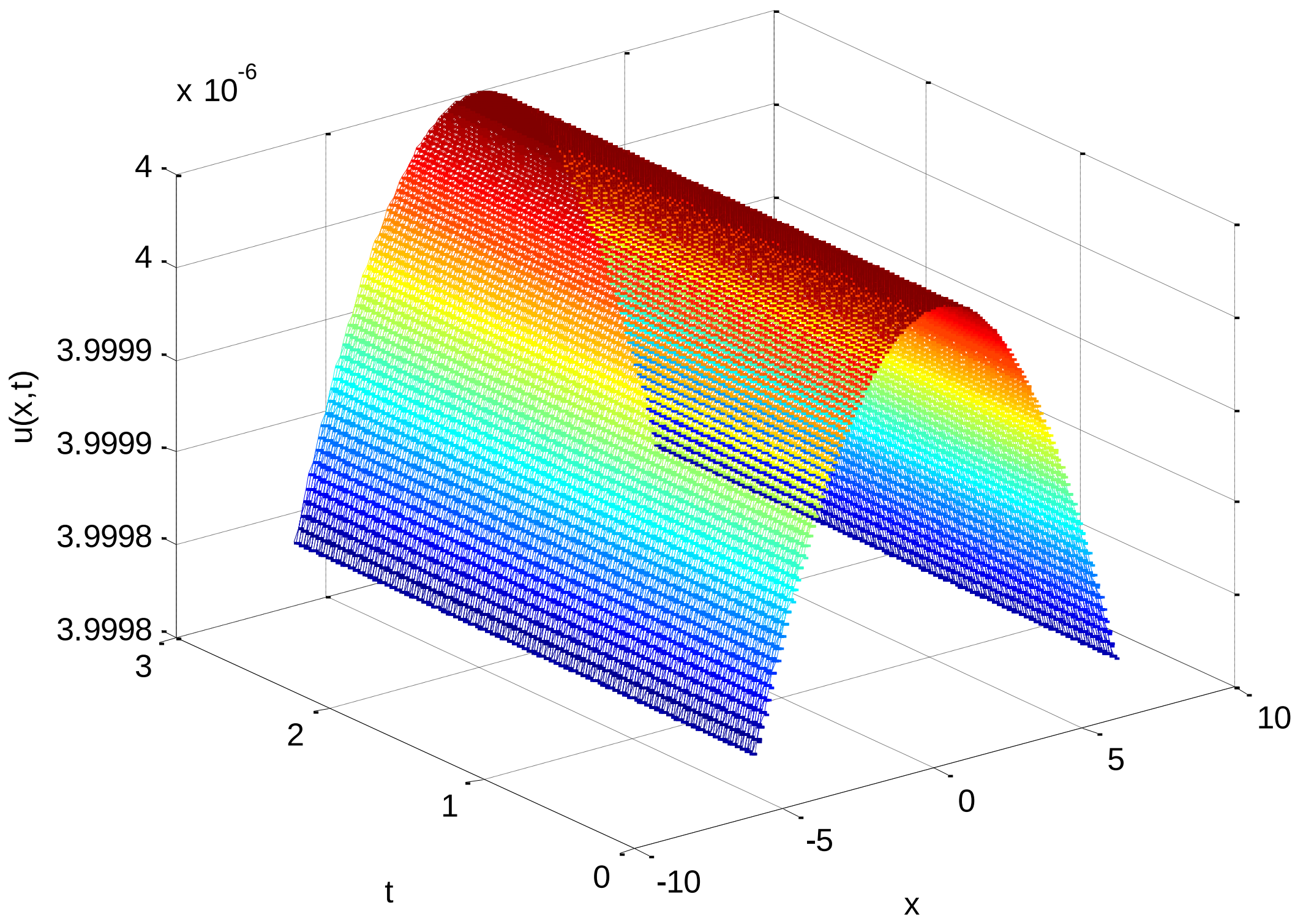

**Figure 1: Mesh plot for numerical solution (Lax), k=0.001**

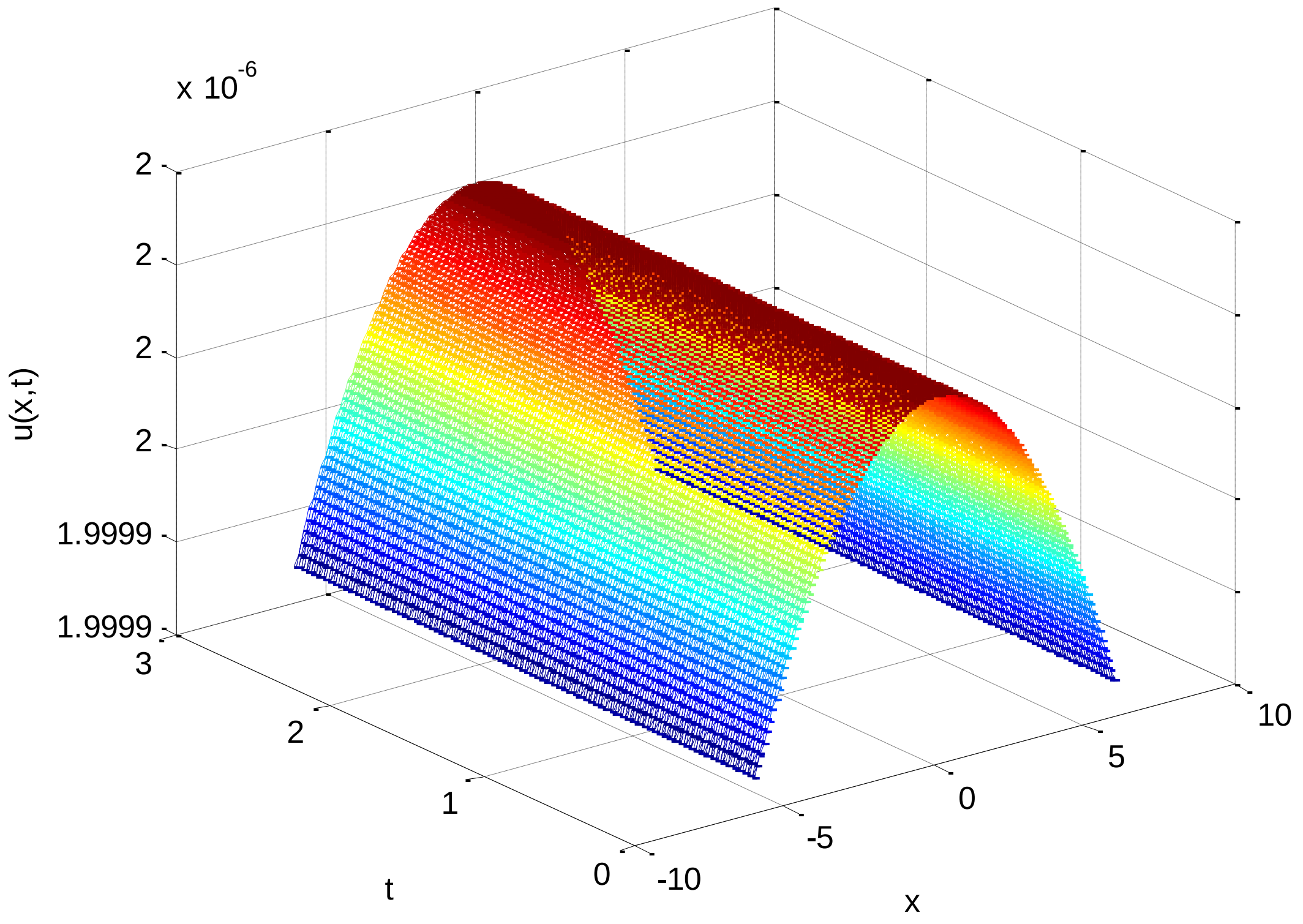

**Figure 2: Mesh plot for numerical solution (SK), k=0.001**

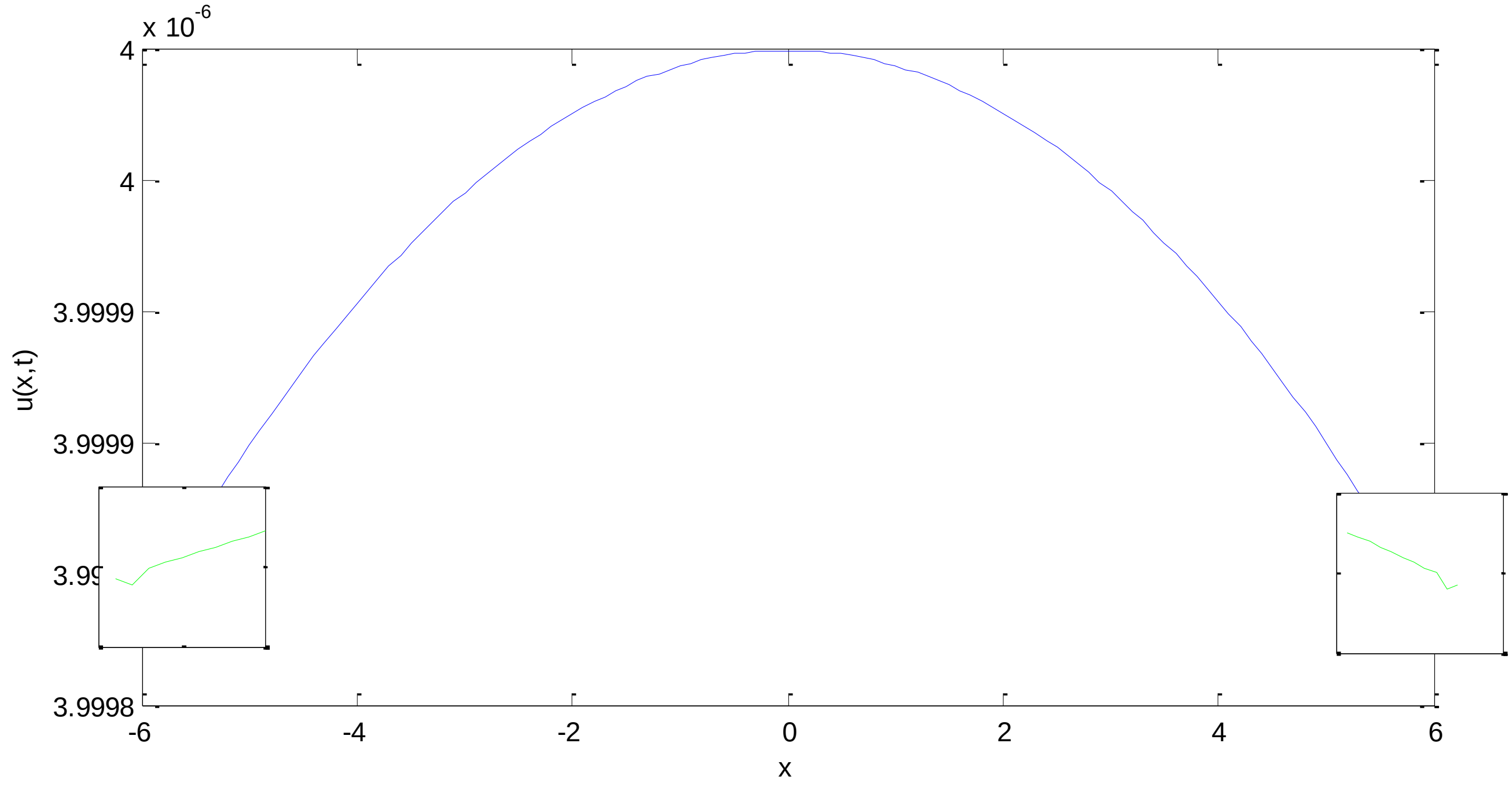

**Figure 3: Error distribution for the solution of Lax for GA and non optimal**

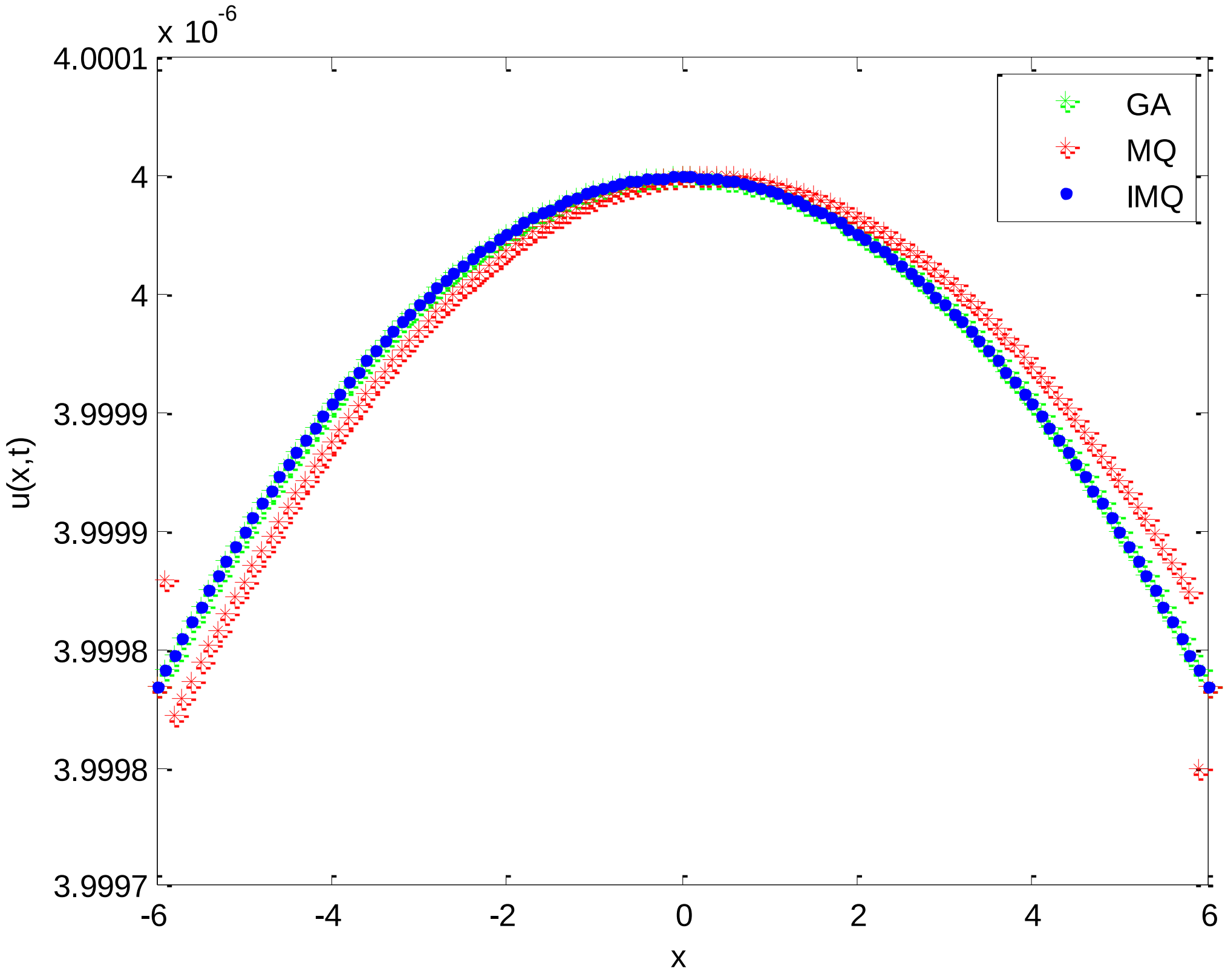

**Figure 4: Solution of Lax for *t=2, k=0.001***

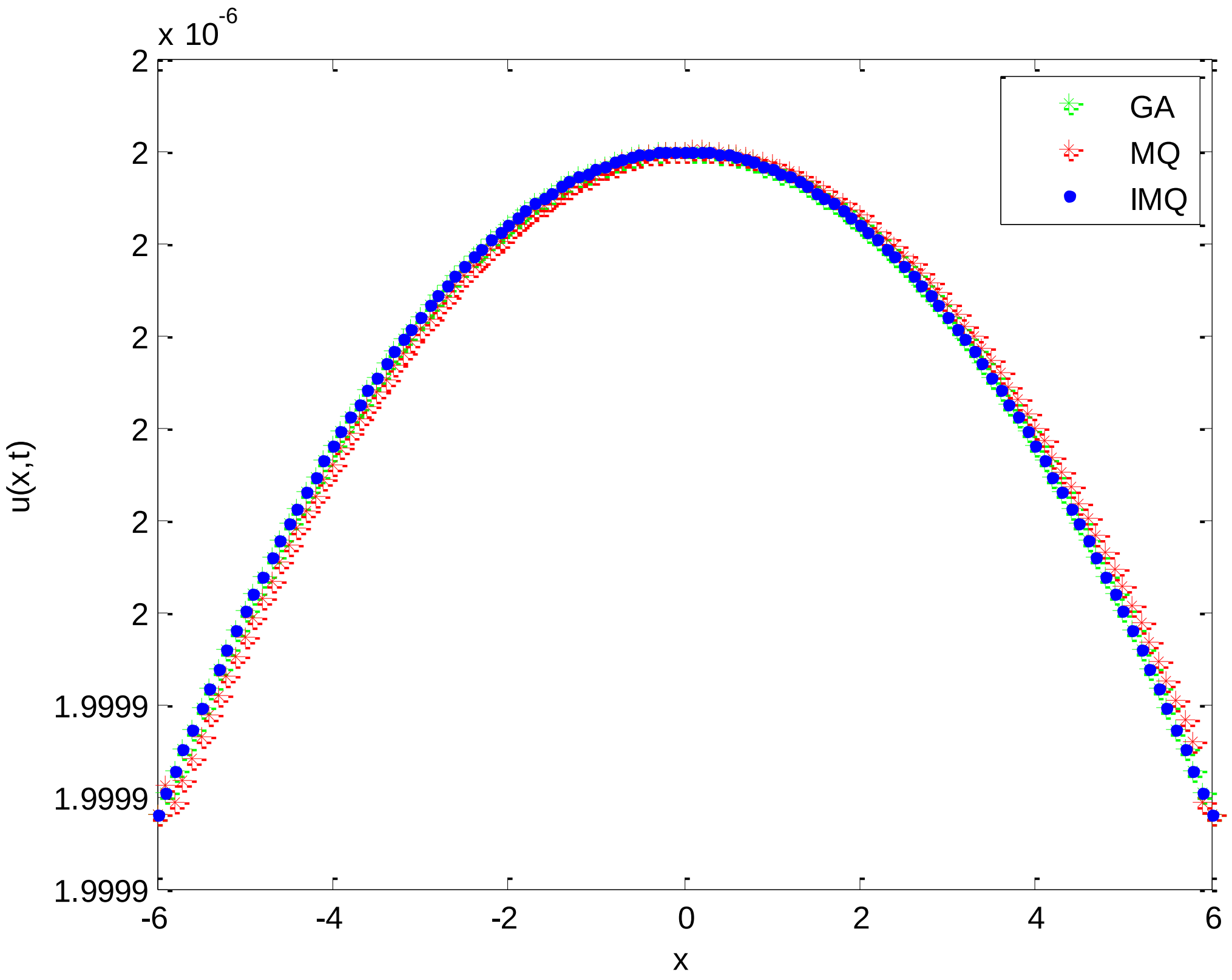

**Figure 5: Solution of SK for *t=2, k=0.001***

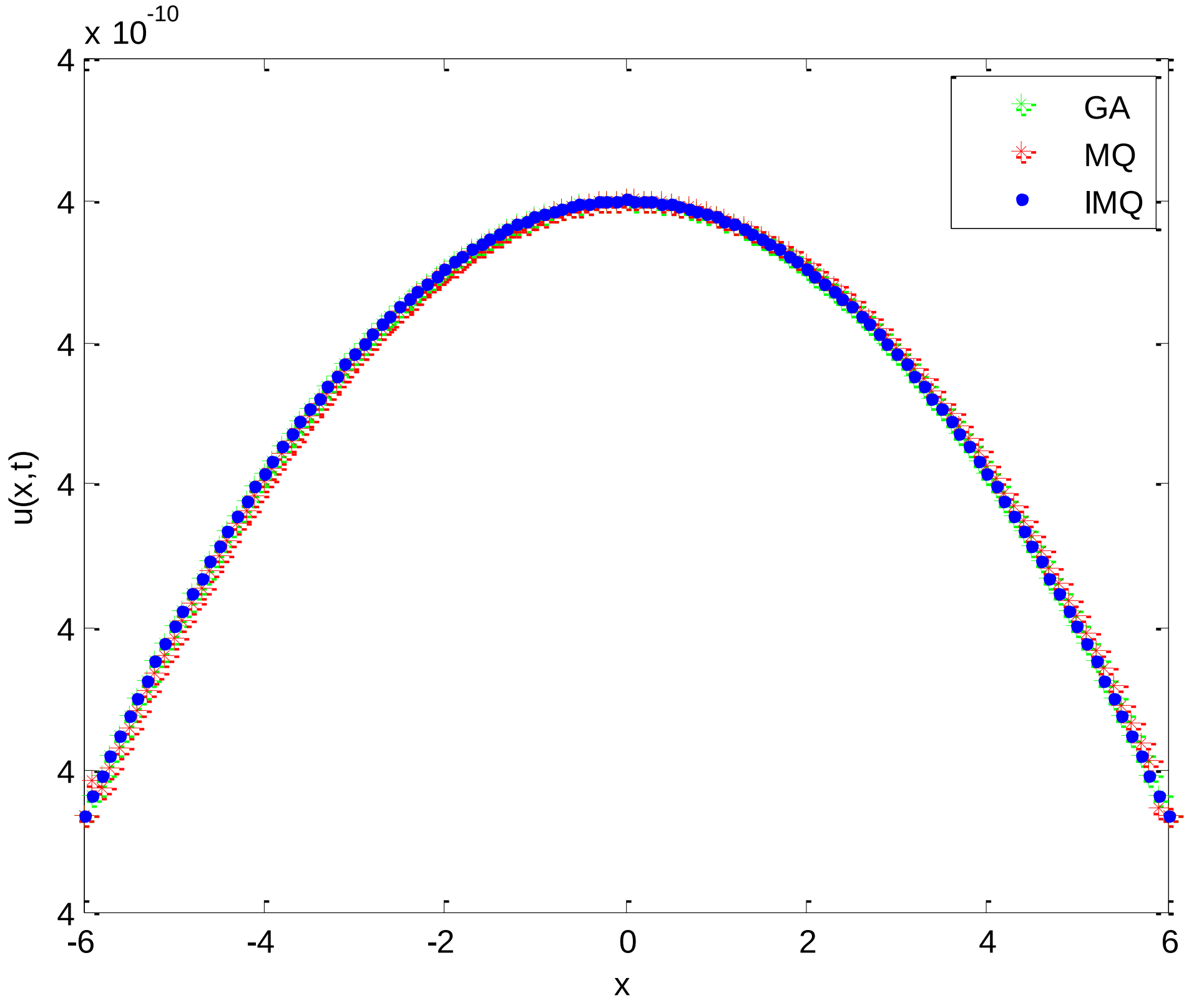

**Figure 6: Solution of Lax for *t=2, k=0.00001***

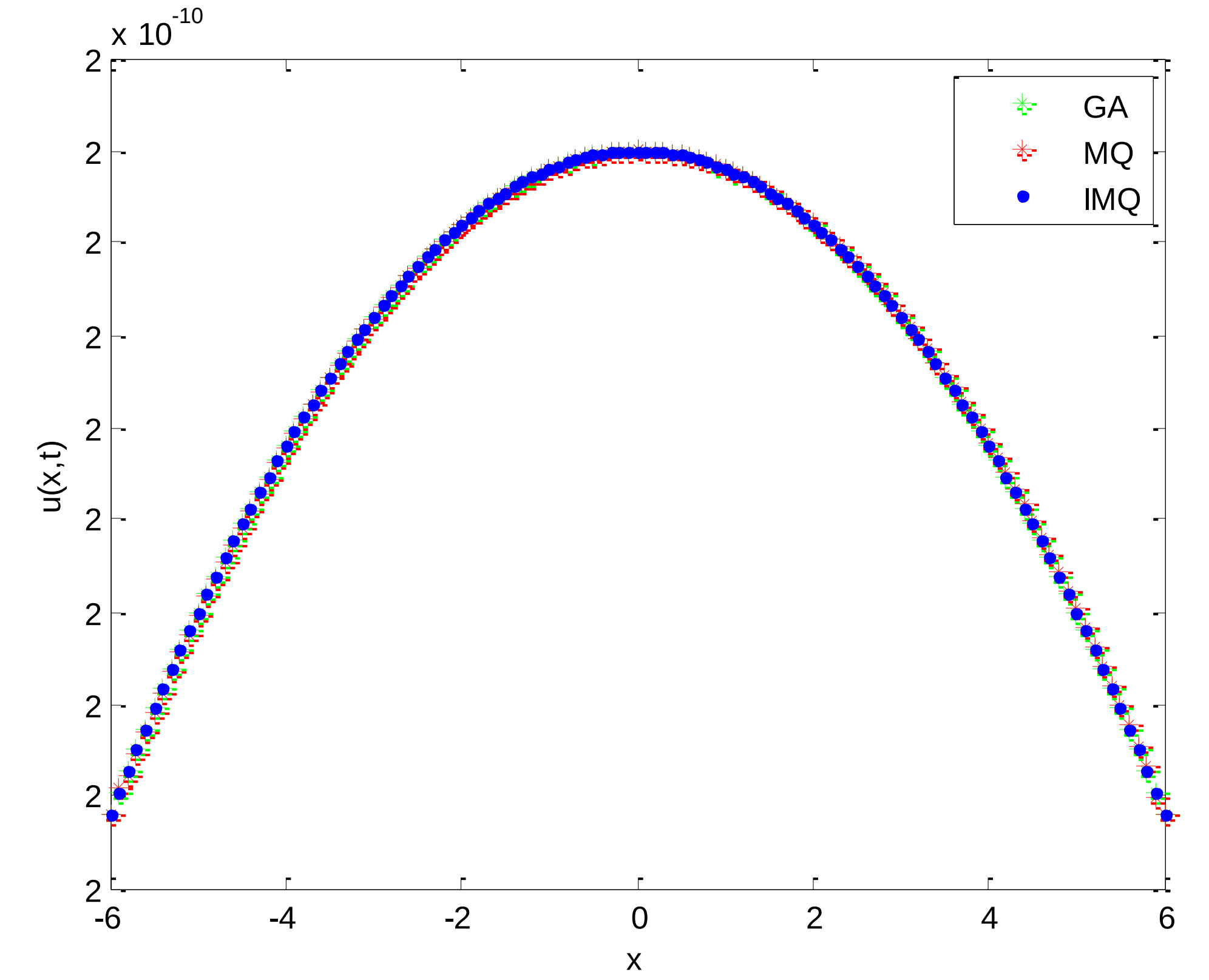

**Figure 7: Solution of SK for *t=2, k=0.00001***

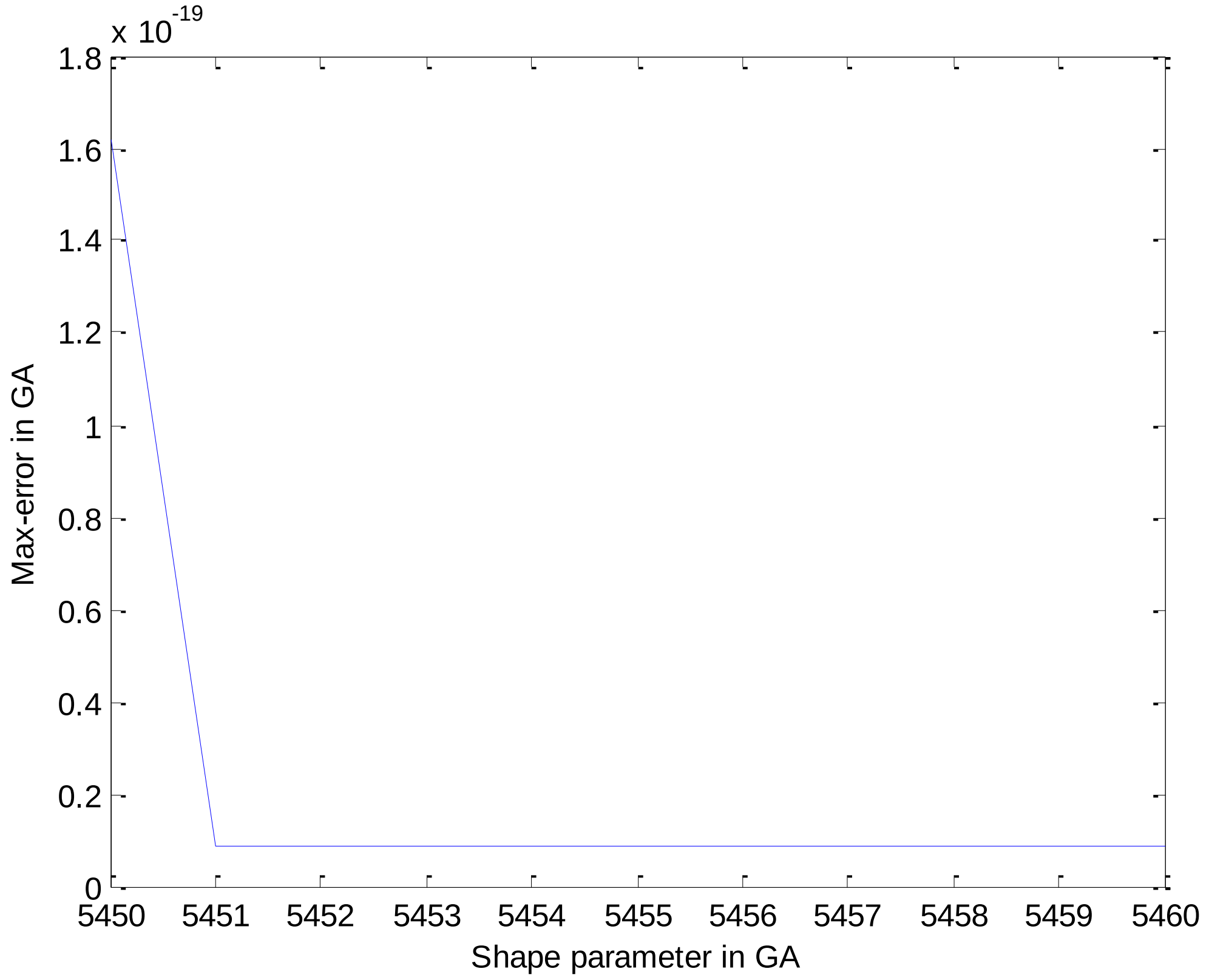

**Figure 8: Shape vs. Max Error, Lax, *k=0.001***

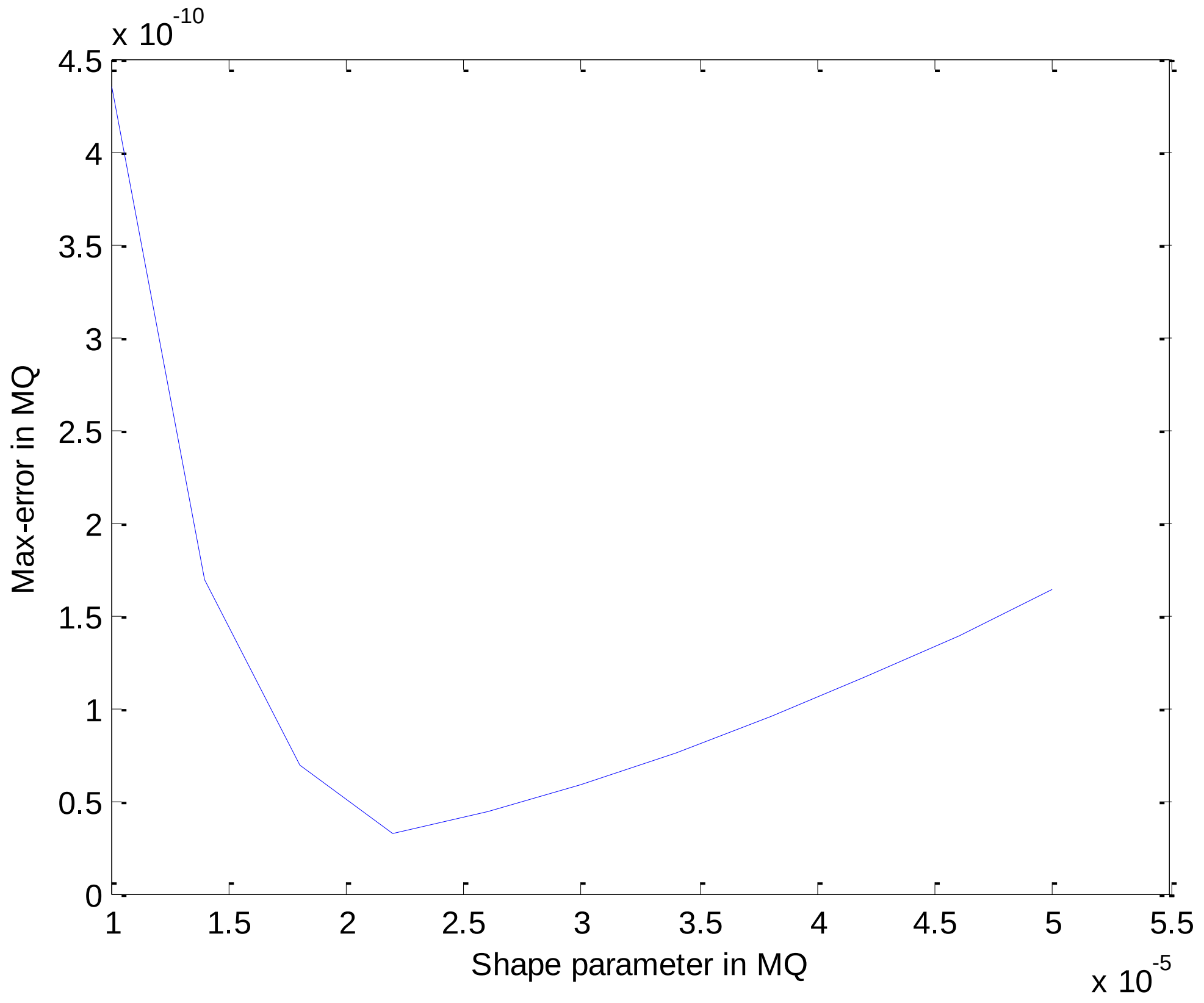

**Figure 9: Shape vs. Max Error, Lax, *k=0.001***

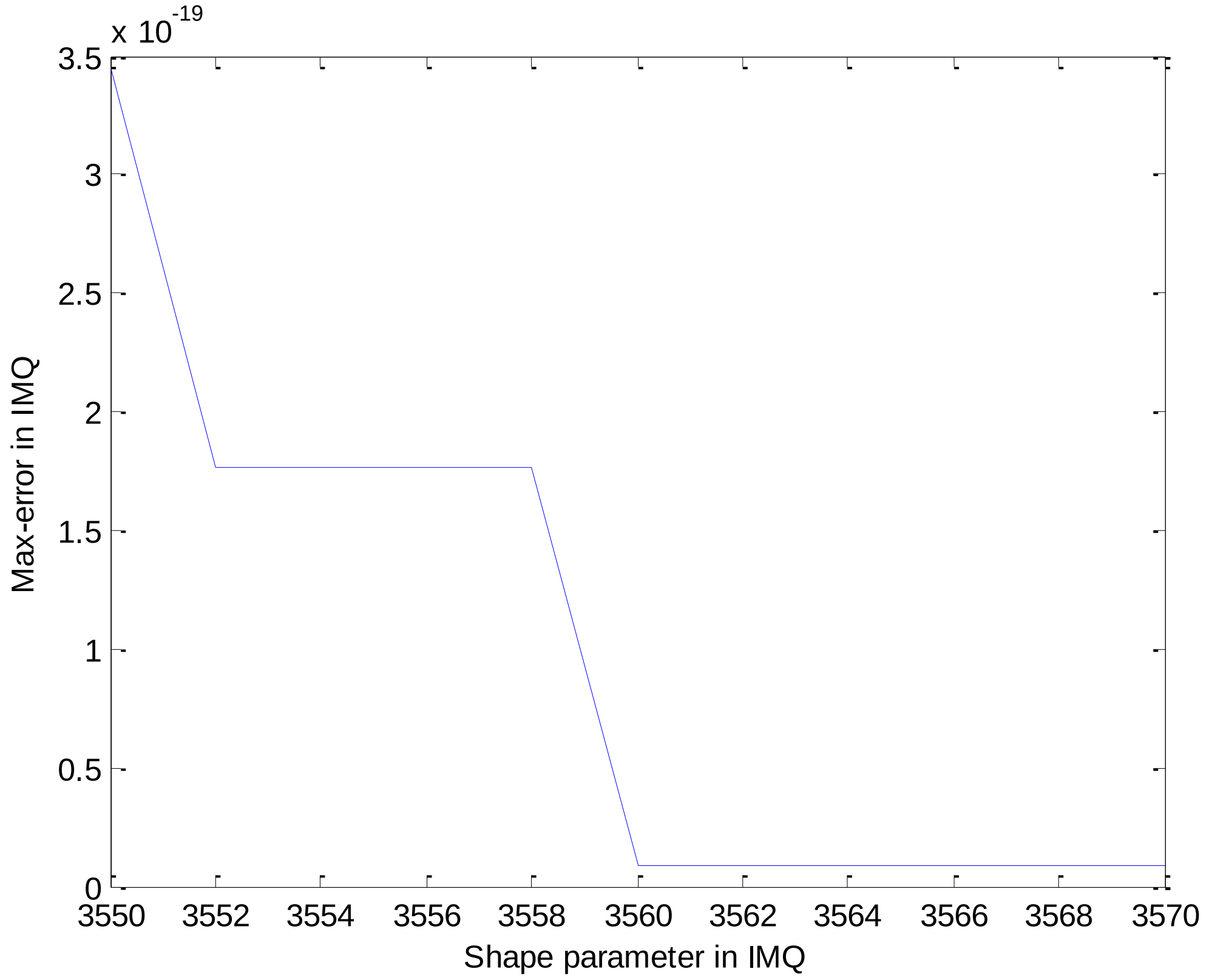

**Figure 10: Shape vs. Max Error, Lax, *k=0.001***

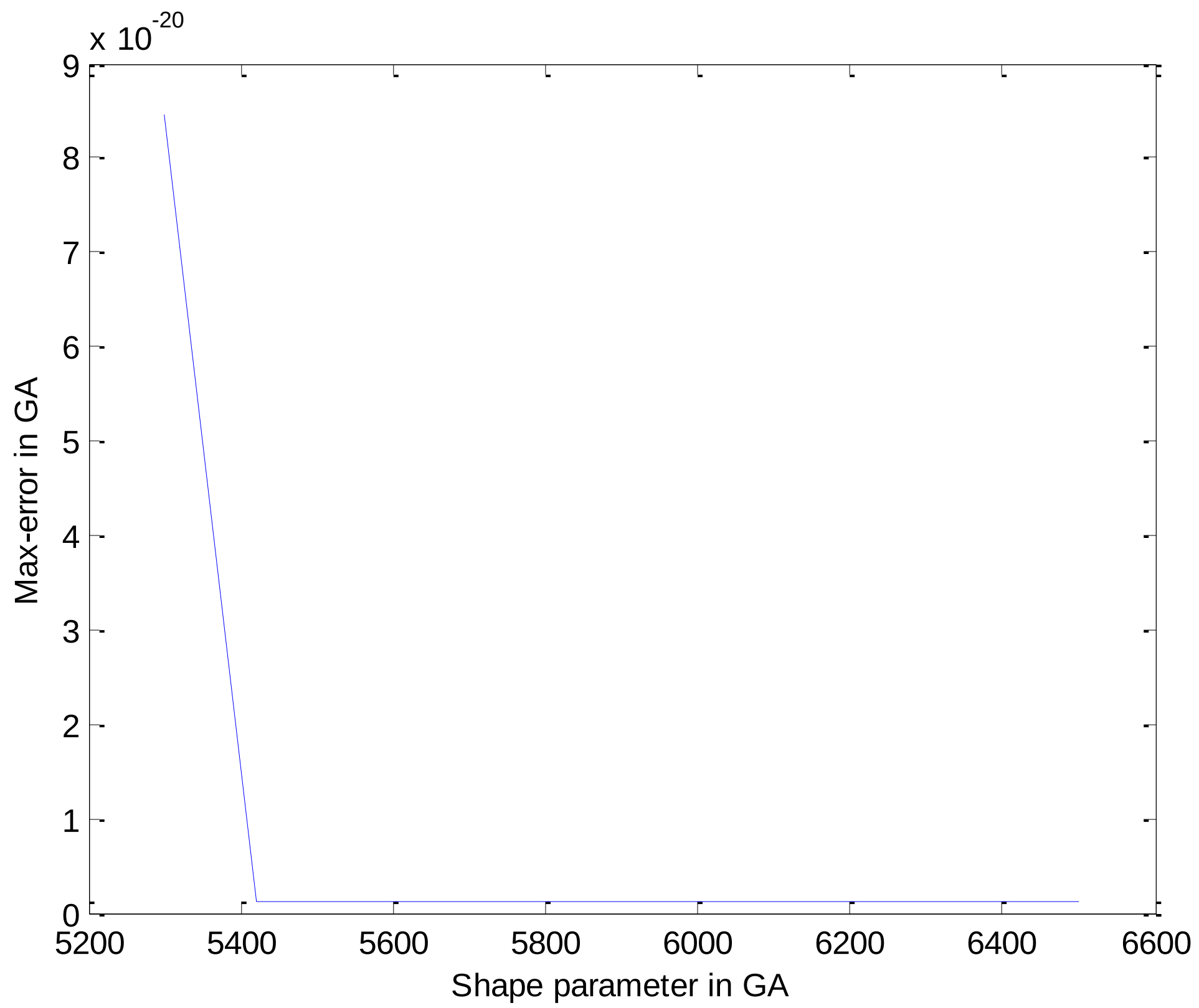

**Figure 11: Shape vs. Max Error, SK, *k=0.001***

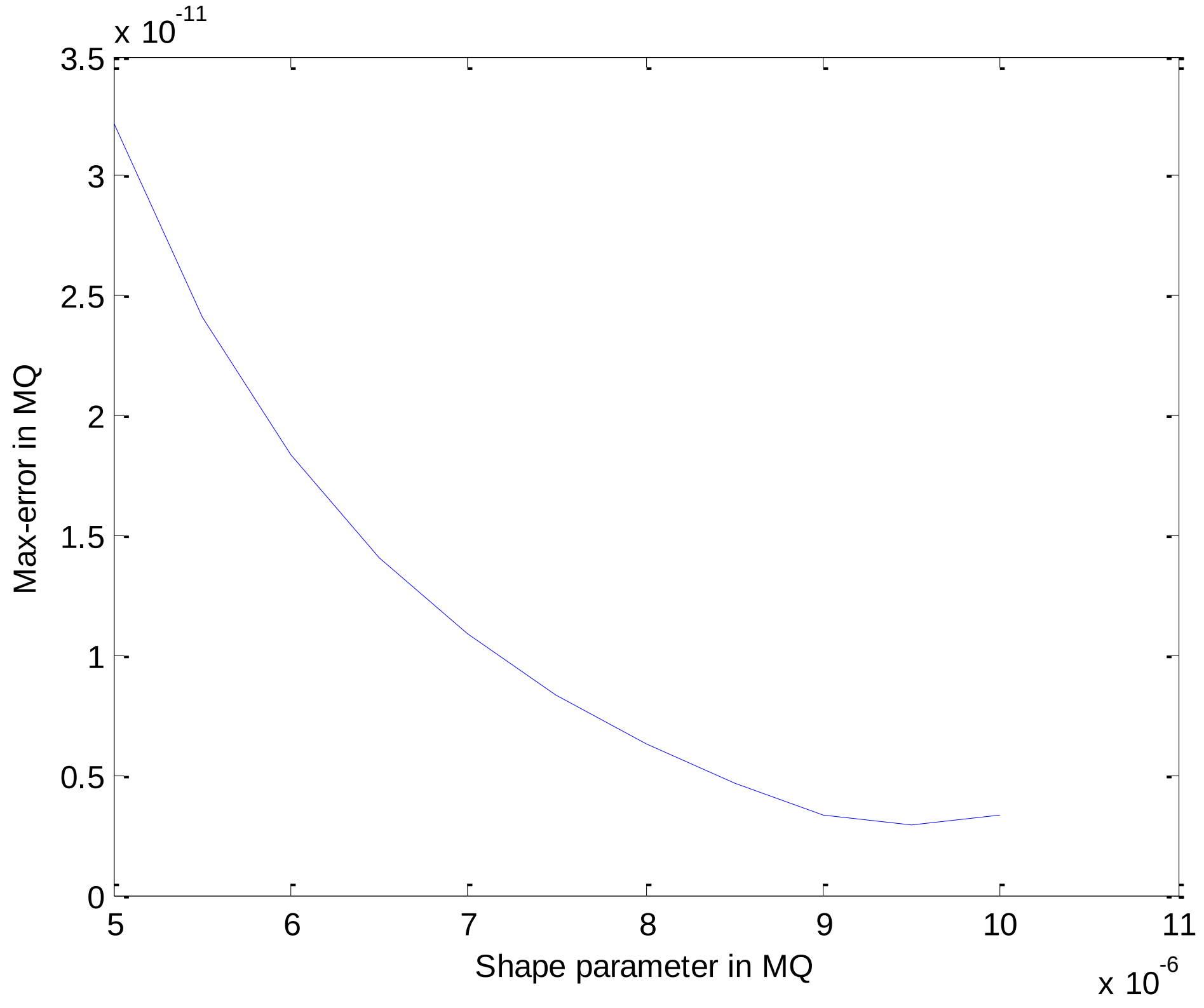

**Figure 12: Shape vs. Max Error, SK, *k=0.001***

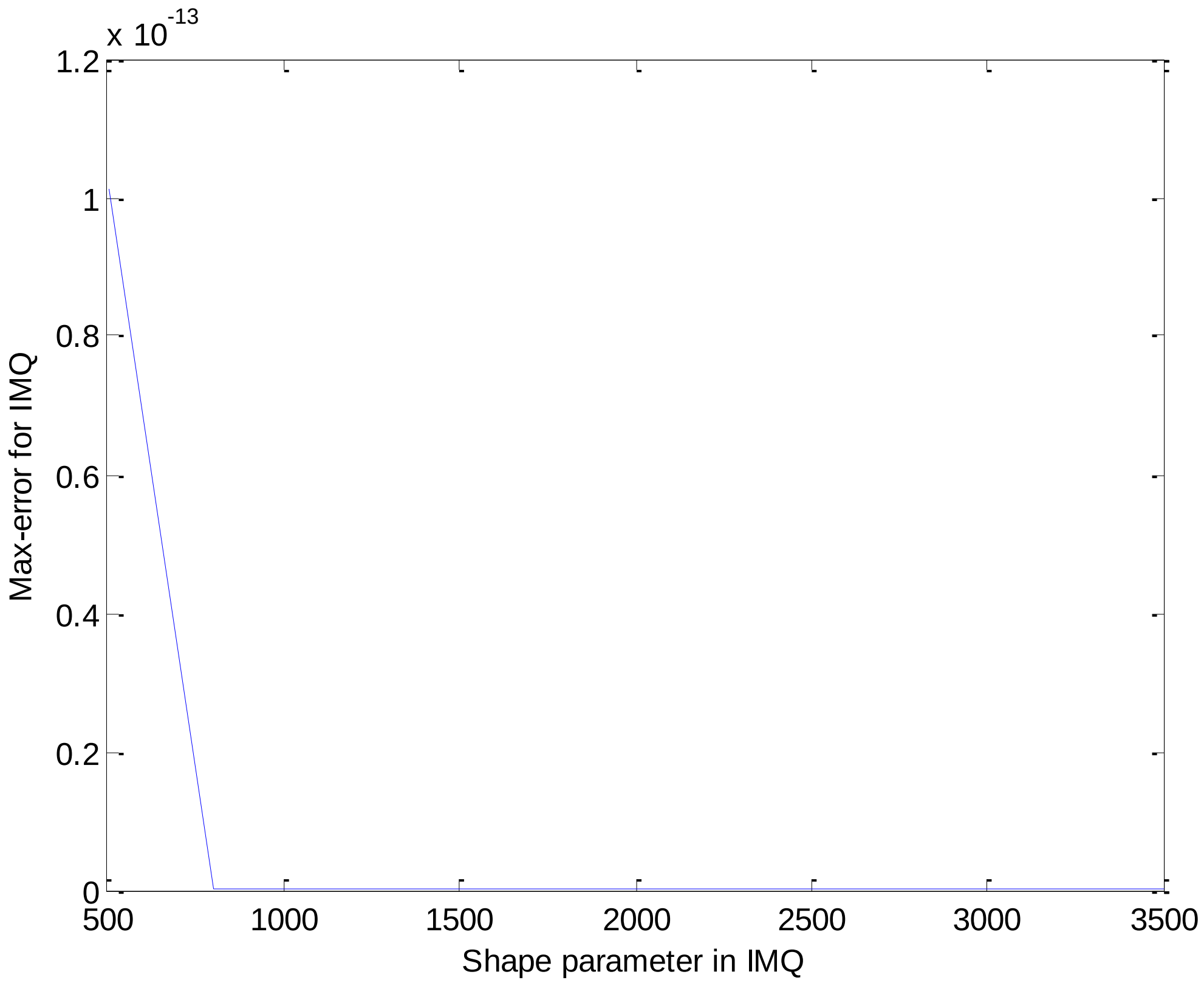

**Figure 13: Shape vs. Max Error, SK, *k=0.001***

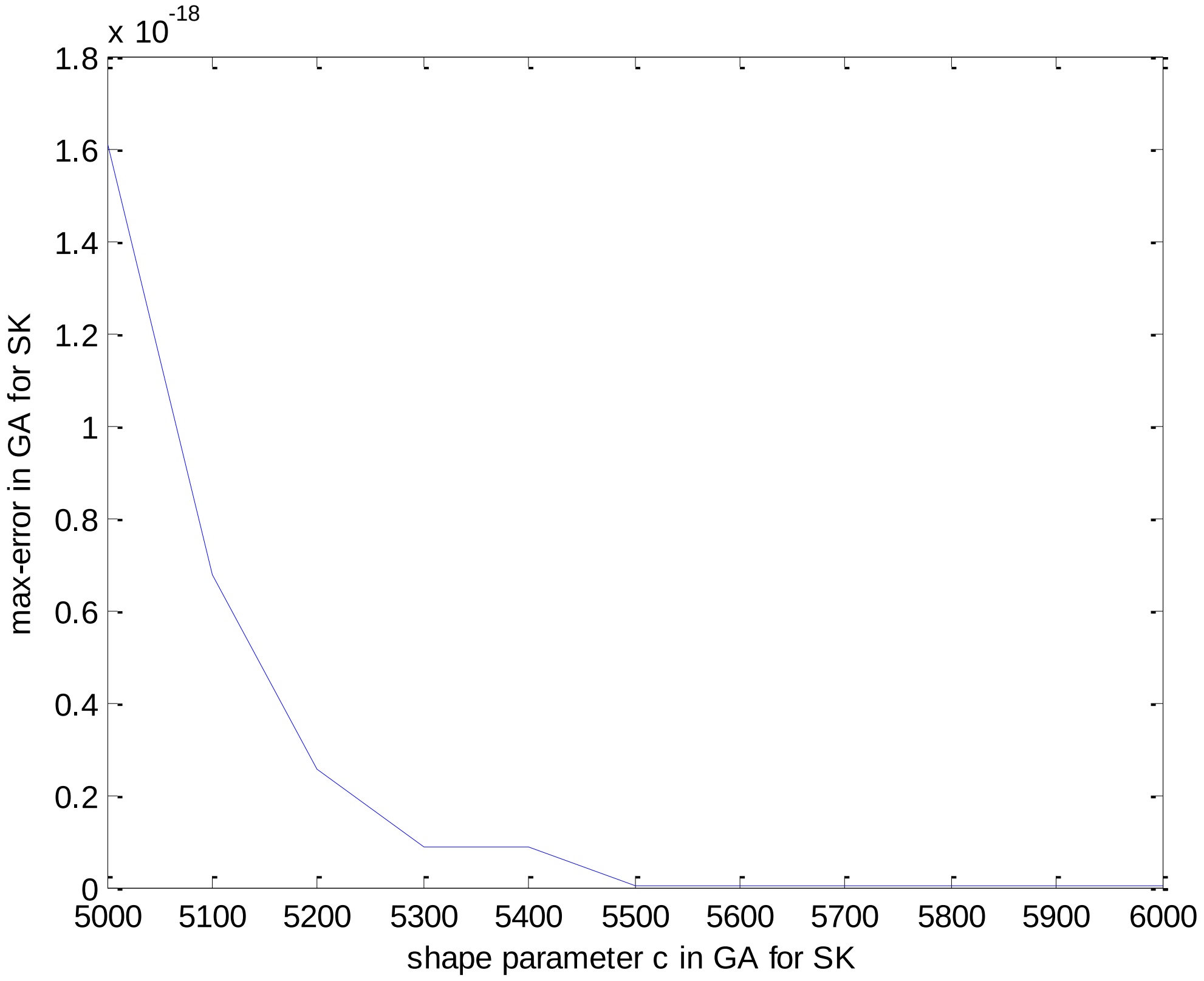

**Figure 14: Shape vs. Max Error, SK, *k=0.001***

**Table 1: Error Comparison for gfKdV (Lax) with $x_0 = 0, k = 0.001$**

| | Shape Parameter | RBF | Max-error | $L_2$-error | RMS error | Cond No. |
|---|---|---|---|---|---|---|
| LAX | 5451 | GA | $8.4703\times10^{-21}$ | $1.5962\times10^{-20}$ | $4.5886\times10^{-21}$ | 1 |
| | $2\times10^{-5}$ | MQ | $4.1421\times10^{-11}$ | $5.0851\times10^{-11}$ | $1.4619\times10^{-11}$ | $1.0174\times10^{4}$ |
| | 3987 | IMQ | $8.4703\times10^{-21}$ | $1.5962\times10^{-20}$ | $4.5886\times10^{-21}$ | $1.2385\times10^{17}$ |
| SK | 5451 | GA | $1.2705\times10^{-21}$ | $1.9774\times10^{-21}$ | $5.6847\times10^{-22}$ | 1 |
| | $10^{-5}$ | MQ | $3.2896\times10^{-12}$ | $5.9692\times10^{-12}$ | $1.7160\times10^{-12}$ | $1.0174\times10^{4}$ |
| | 3987 | IMQ | $1.2705\times10^{-21}$ | $1.9774\times10^{-21}$ | $5.6847\times10^{-22}$ | $4.0915\times10^{17}$ |

**Table 2: Error Comparison for gfKdV (Lax) with $x_0 = 0, k = 0.00001$**

| | Shape Parameter | RBF | Max-error | $L_2$-error | RMS error | Condition No. |
|---|---|---|---|---|---|---|
| LAX | 4790 | GA | $1.0340\times10^{-23}$ | $3.2697\times10^{-24}$ | $9.3998\times10^{-25}$ | 1 |
| | $8.4\times10^{-8}$ | MQ | $4.6374\times10^{-20}$ | $8.6520\times10^{-20}$ | $2.4873\times10^{-20}$ | $1.0170\times10^{4}$ |
| | 1432 | IMQ | $1.0340\times10^{-23}$ | $4.6241\times10^{-24}$ | $1.3293\times10^{-24}$ | $2.3602\times10^{17}$ |
| | $3\times10^{-7}$ | IMQ | $2.2618\times10^{-10}$ | $3.3818\times10^{-10}$ | $9.7220\times10^{-11}$ | 1 |
| SK | 4792 | GA | $5.1699\times10^{-24}$ | $1.6349\times10^{-24}$ | $4.6999\times10^{-25}$ | 1 |
| | $3.66\times10^{-8}$ | MQ | $4.3996\times10^{-21}$ | $8.2079\times10^{-21}$ | $2.3596\times10^{-21}$ | $1.0170\times10^{4}$ |
| | 1324 | IMQ | $5.1699\times10^{-24}$ | $1.2664\times10^{-23}$ | $3.6405\times10^{-24}$ | $1.7591\times10^{18}$ |
| | $25\times10^{-8}$ | IMQ | $5.2398\times10^{-11}$ | $1.1609\times10^{-10}$ | $3.3373\times10^{-11}$ | 1 |

**Table 3: Table of Conserved Quantities for gfKdV (Lax)**

| RBF | $I_1$ | $I_2$ |
|---|---|---|
| GA | $4.7999\times10^{-5}$ | $2.5812\times10^{-16}$ |
| MQ | $4.7999\times10^{-5}$ | $2.5812\times10^{-16}$ |
| IMQ | $4.7999\times10^{-5}$ | $2.5812\times10^{-16}$ |

**Table 3: Table of Conserved Quantities for gfKdV (SK)**

| RBF | $I_1$ | $I_2$ |
|---|---|---|
| GA | $2.4000\times10^{-5}$ | $3.2263\times10^{-17}$ |
| MQ | $2.4000\times10^{-5}$ | $3.2263\times10^{-17}$ |
| IMQ | $2.4000\times10^{-5}$ | $3.2263\times10^{-17}$ |